\documentclass[preprint, 12pt]{elsarticle}
\usepackage{url} 
\usepackage{soul}
\usepackage{amsmath}
\usepackage{amssymb}\def\d{{\, \rm d}}

\begin{document}

\begin{frontmatter}

\title{Probabilistic Eddy Identification with Uncertainty Quantification}

%
%

\author[label1]{Jeffrey Covington}
\ead{Jeffrey.Covington@nau.edu}
\affiliation[label1]{organization={Department of Mathematics, University of Wisconsin-Madison},
            addressline={480 Lincoln Dr.},
            city={Madison},
            postcode={53706},
            state={WI},
            country={United States of America}}
            
\author[label1]{Nan Chen \corref{cor1}}
\ead{chennan@math.wisc.edu}

\author[label2,label3]{Stephen Wiggins}
\ead{wiggins@usna.edu}
\affiliation[label2]{organization={School of Mathematics, University of Bristol},
            addressline={Fry Building, Woodland Road},
            city={Bristol},
            postcode={BS8 1UG},
            country={United Kingdom}}           
            
\author[label3]{Evelyn Lunasin}
\ead{lunasin@usna.edu}
\affiliation[label3]{organization={Department of Mathematics, United States Naval Academy},
            addressline={Chauvenet Hall, 572C Holloway Road},
            city={Annapolis},
            postcode={21402-5002},
            state={MD},
            country={United States of America}}

\begin{abstract}
Mesoscale eddies are critical in ocean circulation and the global climate system. Standard eddy identification methods are usually based on deterministic optimal point estimates of the ocean flow field. However, uncertainty exists in estimating the flow field due to noisy, sparse, and indirect observations and turbulent flow models. Because of the intrinsic strong nonlinearity in the eddy identification diagnostics, even a small uncertainty in estimating the flow field can cause a significant error in the identified eddies. This paper presents a general probabilistic eddy identification framework that adapts existing identification methods to incorporate uncertainty into the diagnostic, emphasizing the interaction between the uncertainty in state estimation and the nonlinearity in diagnostics for affecting the identification results. The probabilistic eddy identification framework starts by sampling an ensemble of flow realizations from the probabilistic state estimation, followed by applying traditional nonlinear eddy diagnostics to individual realizations. The corresponding eddy statistics are then aggregated from the diagnostic results based on these realizations. The framework is applied to a scenario mimicking the Beaufort Gyre marginal ice zone, where large uncertainty appears in estimating the ocean field using Lagrangian data assimilation with sparse ice floe trajectories. The skills in counting the number of eddies and computing the probability of each eddy event are significantly improved under the probabilistic framework. Notably, incorporating the nonlinear propagation of uncertainty in diagnostics provides a more accurate mean estimate than standard deterministic methods in estimating eddy lifetime. It also facilitates uncertainty quantification in inferring such a crucial dynamical quantity.
\end{abstract}


\end{frontmatter}

%
%

%


%
%
%
%

\section{Introduction}

Oceanic eddies are dynamic rotating structures in the ocean. These eddies move and evolve over time and form at many spatial and temporal scales. Mesoscale ocean processes, on the order of 50 km to 500 km, are a vital component of the ocean \cite{morrow2012recent, robinson2012eddies}. Mesoscale eddies are the major drivers of the transport of momentum, heat, and mass \cite{zhang2014oceanic, chelton2007global, cheng2014statistical, yankovsky2022influences} as well as biochemical and biomass transport and production \cite{sandulescu2007plankton,angel1983eddies, prairie2012biophysical, gaube2013satellite, mcgillicuddy1998influence} in the ocean. The study of ocean eddies is increasingly important as eddies play a vital role in the rapidly changing polar regions and the global climate system \cite{beech2022long}.

Eddies can be identified, tracked, and studied through both in-situ and remote observations. Satellite altimetry, which measures sea surface height (SSH), is the most widely used source for global and regional eddy datasets \cite{fu2010eddy, morrow2012recent, le2001ocean, pegliasco2021meta3,abdalla2021altimetry}. Despite being a major revolution in tracking eddies, satellite measurements may have relatively low spatial and temporal resolution and suffer from altimetric noise \cite{larnicol1995mean}. In the polar regions, for example, the presence of sea ice dramatically limits the usefulness of altimetry data for eddy identification \cite{von2022eddies,manley1985mesoscale, kubryakov2021large, kozlov2019eddies}. Many other sources of data in the Arctic, including in-situ observations \cite{timmermans2008eddies, zhao2016evolution} and satellite-derived observations of sea ice motion \cite{lopez2019ice, lopez2021library, manucharyan2022spinning}, are used as alternatives for inferring the eddies.

Due to the complex spatial and dynamical structure of eddies, there is no universal eddy identification criterion. Instead, many eddy diagnostics have been developed to characterize and automatically identify eddies from different types of observational data for use in various applications \cite{d2009comparison, souza2011comparison}. The Okubo-Weiss (OW) parameter is a widely used approach based on the physical properties of the ocean flow \cite{okubo1970horizontal, weiss1991dynamics}. Another set of widely used approaches is based on the geometry of the ocean velocity field \cite{nencioli2010vector} or on contours of the measured SSH \cite{chelton2011global, faghmous2015daily}, which approximately correspond with streamlines of geostrophic ocean flow. Wavelet analysis applied to SSH measurements has also been used for eddy identification from altimetry data \cite{doglioli2007tracking}. Such approaches have the advantage of being directly applicable to altimetry data. On the other hand, there are Lagrangian approaches to eddy identification \cite{branicki2011lagrangian}, such as Finite-time Lyapunov exponents \cite{nolan2020finite} and Lagrangian descriptors \cite{mancho2013lagrangian, vortmeyer2016detecting, vortmeyer2019comparing} which incorporate dynamical information from the ocean. These approaches reflect the dynamic rather than static structure of eddies in the ocean flow. Machine learning has also been applied to eddy identification \cite{franz2018ocean, xu2019oceanic, wang2022data}. The choice of which eddy diagnostics to utilize depends on the intended use case, the available observational data, and other practical considerations.

One major challenge in eddy identification is that the estimated states of the flow field often have inherent uncertainty due to the turbulent flow dynamics, the measurement noise, indirect measurements, and the limited spatiotemporal resolutions in observations \cite{majda2018model, mignolet2008stochastic, palmer2001nonlinear}. Because of the intrinsic strong nonlinearity in the eddy identification diagnostics, even a small uncertainty in estimating the flow field can cause a significant error in the identified eddies \cite{resseguier2017geophysical1, resseguier2017geophysical2, resseguier2017geophysical3}. Altimetry data, for example, does not directly measure the velocity field, which can pose problems for some eddy diagnostics that are sensitive to the altimetric noise \cite{chelton2011global, d2009comparison, souza2011comparison}. Another example is using sea ice floe trajectories as a source of Lagrangian data when altimetry data is unreliable \cite{lopez2019ice, lopez2021library, manucharyan2022spinning}, where uncertainties due to sparse observations can affect the identification skill. 

Data assimilation, which combines noisy, sparse, and indirect observations with a dynamical model, is often utilized to improve the state estimation of the entire ocean as a prerequisite of eddy identification. Particularly, data assimilation incorporates the sparse observations into the spatiotemporal dependence of the ocean field provided by the dynamical model to recover the flow field in the regions not covered by observations. This offers the recovery of the entire spatiotemporal ocean field that facilitates the study of eddies. Notably, the optimal point estimate of the ocean state from data assimilation is usually given by the posterior mean value at each spatiotemporal grid point. It is commonly used as the recovered ocean field in the reanalysis product, which is then naturally adopted for eddy identification \cite{neveu2016historical, song2012application, de2018using, subramanian2013data, xie2020impact}. However, the dynamical and statistical features of the full flow field are characterized not only by the point estimates but also by the associated uncertainty in the state estimation. In the presence of uncertainty, the reconstructed spatiotemporal field from even the optimal point estimates may become insufficient to capture the crucial properties of the underlying ocean dynamics. Large biases appear in the subsequent eddy identification since the uncertainty is significantly amplified by the strong nonlinearity in the eddy diagnostics. Therefore, the uncertainty in state estimation of the ocean needs to be explicitly incorporated into the method of eddy identification. The importance of integrating and interpreting uncertainty has been emphasized through visualization \cite{raith2021uncertainty, hollt2014ovis, potter2009ensemble, potter2012quantification}.

This paper presents a probabilistic eddy identification framework that adapts existing identification methods to incorporate uncertainty into the diagnostic. It highlights how uncertainty in state estimation is amplified through the nonlinearity in the diagnostic criteria when identifying the eddies. At a high level, the probabilistic eddy identification framework simulates an ensemble of realizations of the estimated ocean flow field. The difference between these realizations comes from the uncertainty in the state estimation of the ocean, characterized by the probability density function (PDF). Notably, traditional eddy diagnostics can be directly applied to individual realizations of the ocean. The corresponding eddy statistics, such as size and lifetime, can then be aggregated from the diagnostic results based on these realizations. It is worth highlighting that several recently developed reanalysis data sets have included ensemble members as the outcome products \cite{zuo2017new, 81237, o2021cafe60v1, zhou2022development}. They facilitate applying the probabilistic framework to eddy identification and uncertainty quantification.

The probabilistic eddy identification framework has several advantages over purely deterministic methods based on the point estimate of the ocean field. First, by connecting the optimal point estimate at each grid point, the resulting spatiotemporal flow field may lack crucial physics. This is because the optimal point estimate often fails to fully characterize the turbulent fluctuations, the information of which is usually contained in the estimated uncertainty. In contrast, the realizations of the ocean flow sampled from the distribution of the state estimation take into account the uncertainty and are more dynamically consistent with nature. Second, the frequency and intensity of extreme events are often underestimated in the point estimates when the average is taken, leading to biases in computing the eddy statistics. Contrary to the point estimate, the probabilistic framework includes information from the entire distribution of the estimated state. Therefore, it can quantify the uncertainty in state estimation and the associated eddy identification, including assessing the likelihood of extreme events, the confidence intervals on relevant quantities, and a full distribution of the possible spatiotemporal field. Third, the probabilistic framework is flexible, where existing tools for eddy identification can be naturally adapted to such a framework. This can be achieved by using a given eddy diagnostic method for each sampled realization of the ocean field, followed by integrating the results to form a statistical representation of the eddy identification.

The efficacy of the probabilistic eddy identification framework will be demonstrated by applying it to a scenario mimicking the Beaufort Gyre marginal ice zone of the Arctic Ocean. In such a situation, the underlying ocean flow is estimated by observing a limited number of sea ice floe trajectories from the Lagrangian data assimilation, where the uncertainty in the state estimation of the ocean arrives as a natural consequence of the sparse and indirect observations as well as the intrinsic turbulent features of the ocean dynamics \cite{chen2022efficient}. The experiment will illustrate the effectiveness and robustness of the probabilistic eddy identification framework and highlight the difference with its deterministic counterpart.

The remainder of this paper is organized as follows. Section \ref{sec:preliminaries} contains tools for state estimation with uncertainty quantification and eddy diagnostic that are the prerequisites for developing the probabilistic eddy identification method.  Section \ref{sec:framework} describes the probabilistic eddy identification framework. Section \ref{sec:results} demonstrates and evaluates the probabilistic eddy identification framework using numerical experiments, compared with the results based on the deterministic point estimate. The identified eddy statistics, including lifetime and size, in the presence of uncertainty, are also illustrated in this section. The paper is concluded in Section \ref{sec:conclusions}. The Appendix contains technical details of the models and derivations used in the experiments.

\section{State Estimation, Uncertainty Quantification, and Eddy Diagnostics}\label{sec:preliminaries}

This section includes tools for the state estimation of the flow field with uncertainty quantification and eddy diagnostic. These are the prerequisites for developing the probabilistic eddy identification framework shown in Section \ref{sec:framework}. Depending on the available information, different methods can be exploited for state estimation of the flow field. The probabilistic eddy identification framework developed in Section \ref{sec:framework} is adaptive to any state estimation method as long as the method contains a quantification of the uncertainty of the estimated flow state. In the illustration considered in this paper, Lagrangian data assimilation will be utilized for estimating the ocean field.

\subsection{\label{sec:uncertainty} Uncertainty in state estimation}

Uncertainty is ubiquitous in practice and affects the state estimation of the flow field. Appropriate quantification of the uncertainty is the prerequisite for skillful eddy identification.

Many sources of uncertainty contribute to the state estimation of the flow field. First, observational noise, such as altimetric noise in satellite altimetry data, is a common source of uncertainty that directly influences the state estimation of the underlying flow field. Second, the underlying flow field is not directly observed in many practical scenarios. Therefore, assumptions are needed to build the connection between the observational quantities and the state variables of the flow field. Uncertainty is inevitably propagated or amplified through these assumptions. For example, sea surface height measured from satellite altimetry is often utilized to infer the ocean velocity field by assuming the flow is geostrophically balanced \cite{chelton2007global}. However, such a simplification introduces additional uncertainty in characterizing the flow field due to the ignorance of the unbalanced component in the underlying flow field, affecting the eddy diagnostics \cite{chelton2011global}. Third, sparse observations appear in various situations. The sparsity may come from limited spatiotemporal resolutions or due to external factors that intervene in measurements, such as clouds. Notably, Lagrangian observations, such as passive tracers or sea ice floe trajectories, are widely utilized for estimating the flow states \cite{apte2013impact}. As the number of Lagrangian drifters is usually much smaller than the degree of freedom of the flow field, such a scenario also corresponds to the case with sparse observations. Finally, since the flow dynamics are typically multiscale and turbulent, the flow field predicted by a dynamical or stochastic model naturally contains high uncertainty as well.

\subsection{A quick overview of data assimilation}

Data assimilation combines observations with a model to estimate the state of a dynamical system with an appropriate quantification of the uncertainty \cite{kalnay2003atmospheric, lahoz2010data, majda2012filtering, asch2016data}. It is typically based on the Bayesian inference, and the resulting state estimation is given by a distribution called the posterior distribution. On the one hand, observations are noisy and sparse. Inferring the quantity of interest from the observational variables can further amplify the uncertainty. On the other hand, models also have uncertainties due to the assumptions made in model development, the limited spatiotemporal resolution, and the approximations and parameterizations used. Nevertheless, a model incorporating physical knowledge can connect the observational process and the variables of interest. An important feature of data assimilation is how it balances information from the observations and the model. When observations are high quality, data assimilation relies more on these observations and less on the model for state estimation. Conversely, when observations are sparse and noisy, the model plays a more critical role, which helps fill in the informational gap in the observations. Even with sparse observations, they help mitigate any biases or errors introduced by the model. There are a variety of approaches to data assimilation \cite{fukumori2001data, evensen2009data,pham2001stochastic}. Data assimilation has been applied to assist eddy identification \cite{neveu2016historical, song2012application, de2018using, subramanian2013data, xie2020impact}, though the uncertainty was not highlighted in such studies. Lagrangian data assimilation is a special type of data assimilation methods. It exploits the trajectories of moving tracers as observations to recover the underlying flow field \cite{griffa2007lagrangian, blunden2019look, honnorat2009lagrangian, salman2008using, castellari2001prediction}. Previous papers have applied Lagrangian data assimilation to observations of sea ice floe trajectories \cite{chen2021lagrangian, covington2022bridging, chen2022efficient}.

\subsection{\label{sec:da} Lagrangian data assimilation using observed sparse sea ice floe trajectories}

This paper will focus on the eddy identification in a scenario mimicking the marginal ice zone of the Arctic Ocean. Such a study has practical significance and is related to understanding climate change. However, direct observations of ocean field in the marginal ice zone is difficult to obtain due to the interference of sea ice with satellite altimetry of the ocean. Therefore, alternate sources of observational data are needed. A practical approach is to utilize observations of individual ice floe trajectories \cite{lopez2019ice, lopez2021library} to estimate the underlying ocean state and then identify eddies \cite{manucharyan2022spinning}. Yet, due to the inherent uncertainty in the underlying turbulent ocean dynamics and the sparse observations of ice floe trajectories, this scenario provides a challenging test for identifying the eddies in the entire spatiotemporal flow field. To facilitate the task, the state estimation of the underlying ocean flow field is given by the Lagrangian data assimilation results, where the observed sea ice floe trajectories play the role of the Lagrangian drifters. The sparse Lagrangian floe observations are combined with a dynamical model via Lagrangian data assimilation to recover the entire spatiotemporal flow field with uncertainty quantification, which will be used for the probabilistic eddy identification. In such a setup, the uncertainty comes from several sources. First, despite being small, observational error exists in determining the floe locations. Second, the observations are indirect, as the underlying flow field needs to be inferred from the Lagrangian drifters. Third, the sparsity of the observations further introduces uncertainty. Notably, the uncertainty is amplified when the sparse Lagrangian trajectories are utilized to estimate the flow field in the regions with no observations.

Despite the strong nonlinearity in Lagrangian data assimilation, an analytically solvable data assimilation scheme is utilized in this study \cite{chen2018conditional, chen2020efficient}, which provides efficient and effective estimation for state estimation. It avoids many empirical tuning procedures, such as the covariance inflation and localization, in the standard ensemble data assimilation methods. The method, therefore, minimizes potential numerical errors introduced by the data assimilation and facilitates conducting numerical experiments to evaluate the probabilistic eddy identification framework. The details of such a Lagrangian data assimilation are included in appendix \ref{sec:cg}.

\subsection{\label{sec:ow}Eddy criterion using the OW parameter}
The primary purpose of this paper is to demonstrate the probabilistic eddy identification framework and provide concrete examples of how uncertainty affects eddy identification. As the goal is not to seek the most appropriate eddy diagnostic baseline approach for specific applications, a simple criterion based on the Okubo-Weiss (OW) parameter \cite{okubo1970horizontal, weiss1991dynamics} is utilized in the following. It is worth highlighting that the probabilistic eddy identification framework developed in this paper is flexible. Any existing eddy diagnostic or criterion applicable to the deterministic eddy identification scenario can be naturally adapted to such a framework.

The OW parameter \cite{okubo1970horizontal, weiss1991dynamics} is a classical and widely-used approach to eddy identification which measures the relative strength of vorticity and strain in the ocean velocity \cite{isern2003identification, isern2006vortices, chelton2007global, cheng2014statistical, petersen2013three}. The OW parameter is derived from the 2-dimensional ocean velocity field $(u_\mathrm{o}, v_\mathrm{o})$ and is defined by
\begin{equation}\label{eqn:OW1}
    \operatorname{OW}(x, y) = s_\mathrm{n}^2 + s_\mathrm{s}^2 - \omega^2
\end{equation}
where the normal strain, the shear strain, and the relative vorticity are given by
\begin{equation}\label{eqn:OW2}
    s_\mathrm{n} = \frac{\partial u_\mathrm{o}}{\partial x} - \frac{\partial v_\mathrm{o}}{\partial y}, \qquad s_\mathrm{s} = \frac{\partial v_\mathrm{o}}{\partial x} + \frac{\partial u_\mathrm{o}}{\partial y}, \qquad\mbox{and} \qquad\omega = \frac{\partial v_\mathrm{o}}{\partial x} - \frac{\partial u_\mathrm{o}}{\partial y}
\end{equation}
respectively. When the OW parameter is negative, the relative vorticity is larger than the strain components, indicating vortical flow. The physical interpretability of the OW parameter is an important advantage of using it in eddy identification.

The eddy criterion used here considers eddies to be regions with a negative OW parameter that is below a certain threshold value. Choosing an appropriate threshold is difficult and, in practice, is region-dependent, but a commonly used threshold is $-0.2 \sigma_\mathrm{OW}$ where $\sigma_\mathrm{OW}$ is the standard deviation of the OW parameter. Eddy cores are considered to be local minima of the OW parameter below the threshold, and eddy boundaries are considered to be closed contours of the OW parameter at the threshold value.

Notably, similar to most eddy diagnostics, the OW parameter is a nonlinear function of the ocean state. Therefore, finding the expected value of the OW parameter based on an ensemble of ocean realizations sampled from the PDF of the estimated state differs from applying the OW parameter to the mean estimate of the ocean. The interaction between nonlinearity and uncertainty can profoundly impact the identification result. The probabilistic approach to eddy identification naturally incorporates the uncertainty in the state estimation and adequately represents the statistical properties of the OW parameter in a coherent and consistent way. The approach has advantages similar to other eddy diagnostics.

\section{\label{sec:framework} The Probabilistic Eddy Identification Framework}

\subsection{Overview}
Eddy identification relies on the estimated ocean state from observational data. In real scenarios, the estimate of the flow field has inherent uncertainty due to measurement noise, observational resolution, model assumptions, and many other factors. Uncertainty can significantly impact eddy identification through the nonlinearity in the eddy diagnostic criteria, so a deterministic point estimate of the ocean may not be sufficient. The probabilistic eddy identification framework developed in this paper incorporates uncertainty into eddy identification. The uncertainty is characterized by an ensemble of realizations sampled from the estimated ocean state. Standard eddy identification methods can be naturally applied to each of these realizations. Then, the statistics of the identified eddies are formed by aggregating the results associated with each.

\subsection{\label{sec:procedure} General procedure of the probabilistic eddy identification framework}
\begin{figure}
    \begin{center}
        \hspace*{-1cm}\includegraphics[width=16cm]{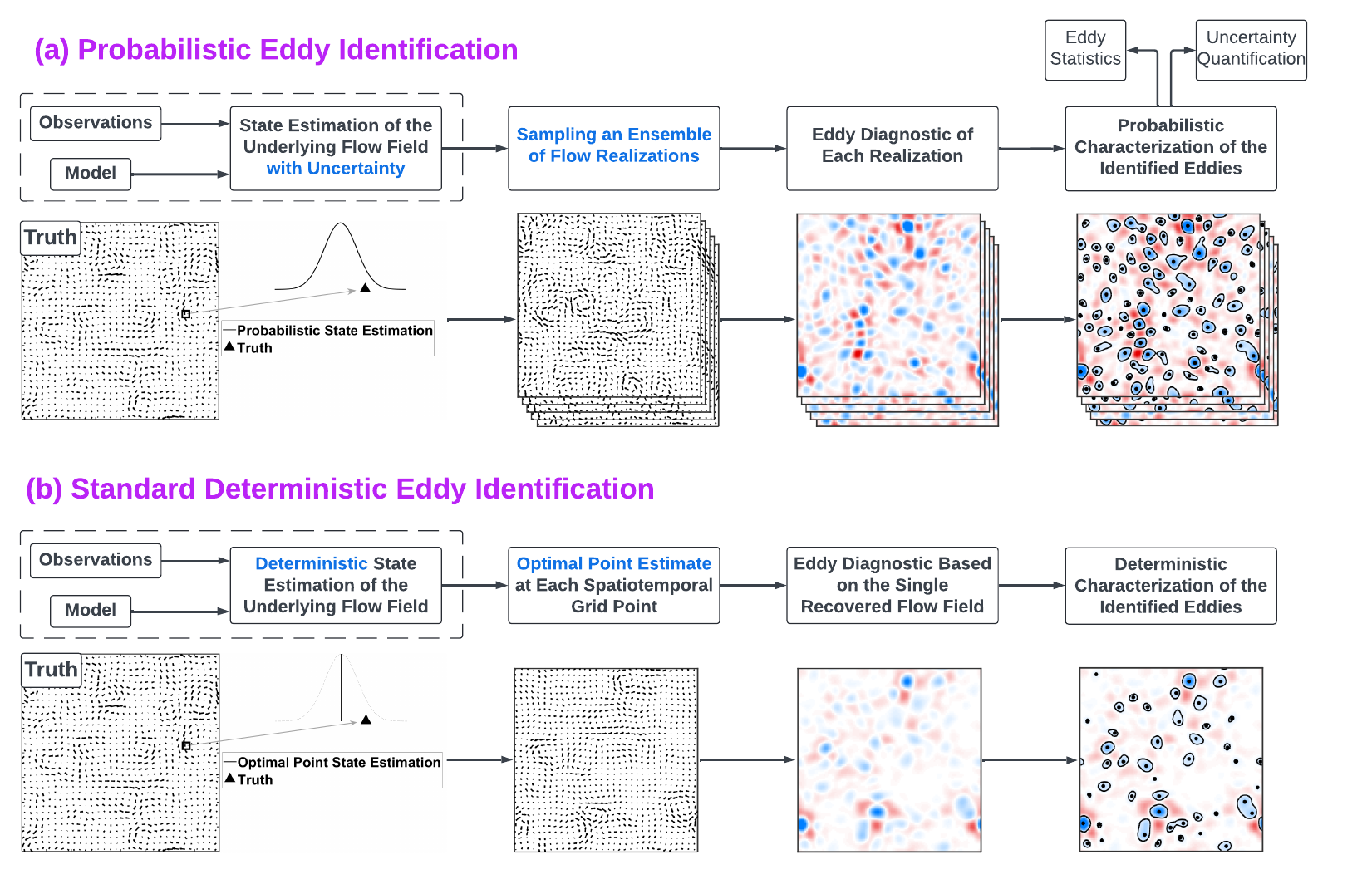}
    \end{center}
    \caption{\label{fig:overview}  Overview of the probabilistic eddy identification framework (Panel (a)), following the general procedure presented in Section \ref{sec:procedure}. The probabilistic eddy identification framework is also compared with the standard deterministic method (Panel (b)).
    }
\end{figure}

The general procedure of the probabilistic eddy identification framework is presented in Panel (a) of Figure \ref{fig:overview}, which contains four key steps.

\noindent\emph{Step 1. State estimation of the flow field with uncertainty quantification.} \\ One prerequisite of eddy identification is recovering the underlying flow field by exploiting the available information from observations and knowledge-based dynamical models. Data assimilation is a natural choice for state estimation to handle uncertainties in observations and models, although other data fusion or inference methods can also be utilized. Notably, the estimated state at each spatiotemporal grid point is given by a PDF, which represents the uncertainty in state estimation and plays a crucial role in characterizing the identified eddies in a probabilistic way.

\noindent\emph{Step 2. Sampling an ensemble of realizations from the distribution of the estimated state.}  \\
Most standard diagnostic approaches are based on a deterministic recovered flow field. However, the recovered flow field from Step 1 is characterized by PDFs. Applying existing realization-based methods directly to PDFs is not straightforward. To facilitate the application of the standard realization-based methods in the presence of uncertainty, an ensemble of realizations of the possible flow field is sampled from these PDFs. Each realization represents one possible scenario of the flow field.

\noindent\emph{Step 3. Applying standard eddy diagnostic to each realization.} \\
Since the PDF of the state estimation is effectively represented by multiple samples, a standard eddy diagnostic method can be naturally applied to each sample, which is a possible realization of the flow field. The diagnostic result remains as an ensemble, with each ensemble member being the eddy diagnostic associated with one sample.

\noindent\emph{Step 4. Aggregating the results associated with each realization to form the probabilistic characterization of the identified eddies.}\\
The eddy diagnostic results from all realizations in Step 3 can then be aggregated to calculate the indicators, including quantifying the resulting uncertainty, for any eddy quantity of interest. This includes the range of the eddy numbers at each time instant, the probability of the occurrence of a specific eddy, the distribution of eddy lifetime, and the uncertainty in calculating eddy size, etc.

\subsection{Comparison between the probabilistic eddy identification framework with the standard deterministic method}
Panel (b) of Figure \ref{fig:overview} shows the general procedure for applying a deterministic eddy identification method. Instead of estimating the flow state with uncertainty quantification, the recovered flow field is given by an optimal point estimate at each spatiotemporal grid point, which can be the mean value computed from the Bayesian inference (e.g., the so-called posterior mean estimate). Then, a standard eddy identification diagnostic is applied to such an estimated flow field.

It is worth noting that the optimal point estimate using the posterior mean estimate usually fails to capture most of the turbulent fluctuations in recovering the underlying flow field. Therefore, it usually fails to capture extreme events and the flow amplitude. As a result, the strength and the number of eddies are underestimated, especially when the uncertainty in the state estimation is significant.

\section{\label{sec:results}Eddy Identification in the Marginal Ice Zone Using Lagrangian Data Assimilation with Ice Floe Trajectories}

The probabilistic eddy identification framework developed in Section \ref{sec:framework} generally applies to a wide variety of eddy identification problems. In this section, the framework is tested on a concrete scenario, mimicking the ocean flow field in the Beaufort Gyre marginal ice zone. Observations of sea ice floe trajectories are combined with a dynamical model to estimate the ocean state and quantify uncertainty. Synthetic data is adopted in this study to validate the method and the results. These synthetic data have been shown in the previous studies to successfully reproduce many observed features of the ice-ocean coupling \cite{covington2022bridging, manucharyan2022spinning, lopez2021library}.

\subsection{\label{sec:setup}Experimental Setup}

The true time evolution of the ocean flow field and the sea ice floe trajectories are generated from an ice-ocean coupled system consisting of a discrete element method (DEM) floe model and a two-layer quasi-geostrophic (QG) ocean model. In the DEM floe model, the individual floe shapes and sizes are drawn from a library of floe observations in the Beaufort Gyre marginal ice zone \cite{lopez2021library}. The domain size is $400$km$\times400$km. In all the experiments except the one in Figure \ref{fig:evofow}, 40 floe trajectories are observed over 100 days. The number of the floes utilized here is analogous to observing floe trajectories in a similarly sized domain in the marginal ice zone over a summer season \cite{lopez2019ice, lopez2021library}. The noise in observing floe locations is $250$m. The details of the model and parameters are included in \ref{sec:ice} and \ref{sec:qg}. By developing an effective stochastic forecast model, an efficient nonlinear Lagrangian data assimilation algorithm is adopted to estimate the state of the ocean. The stochastic forecast model and the Lagrangian data assimilation scheme are shown in \ref{sec:ocean} and \ref{sec:cg}, respectively.

\begin{figure}
    \begin{center}
        \includegraphics[]{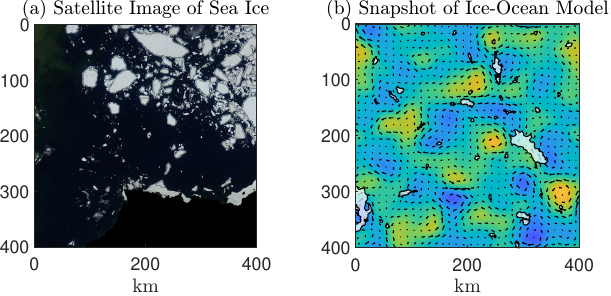}
    \end{center}
    \caption{\label{fig:ice} Panel (a): a processed satellite image of the marginal ice zone in the Beaufort sea in June 2008. Panel (b): a snapshot of the coupled ice-ocean model. The shadings represent the stream function and the arrows show the velocity field of the quasi-geostrophic ocean model.
    }
\end{figure}

Panel (b) of Figure \ref{fig:ice} shows a snapshot of the sea ice model, mimicking the Beaufort Gyre marginal ice zone shown in Panel (a). For simplicity, the floes are assumed to be shape-preserving, where no melting and welding occurs. The experiments here involve a relatively low-concentration region of ice floes. Therefore, the floes are considered to be noninteracting. In addition, atmospheric forcing mainly occurs on a large scale and barely affects the mesoscale motion of the floes, which is often assumed to be known. Since the effect of the atmospheric forcing can roughly be eliminated by applying a Galilean transformation, it is omitted from the following experiment for simplicity. Therefore, the motion of the ice floe is purely driven by the ocean current. While this model simplifies the ice-ocean interactions in some respects, it remains to capture intricate ice-floe dynamics present in real datasets of floe trajectories. These simplifications are made primarily to allow this work to focus on comparing the probabilistic framework with its deterministic counterpart. Note that significant uncertainty in estimating the ocean state still exists due to the use of sparse Lagrangian floe trajectories and a turbulent ocean model.

In all the following experiments, 100 realizations of the ocean field sampled from the state estimation are utilized for computing the statistics under the probabilistic eddy identification.

\subsection{\label{sec:oceanuncertainty}Uncertainty in the state estimation of the ocean flow field}

\begin{figure}
    \begin{center}
        \includegraphics{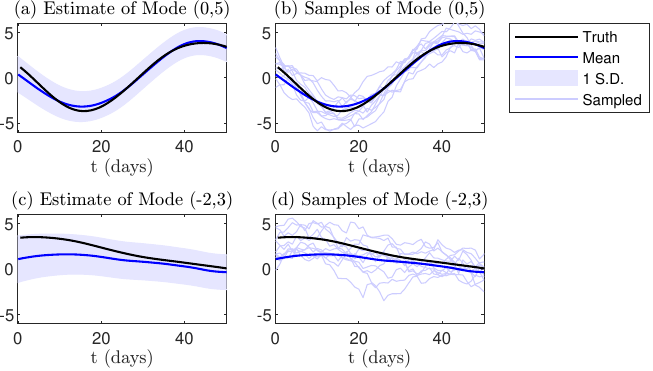}
    \end{center}
    \caption{\label{fig:smoother} The results of the Lagrangian data assimilation for the real parts of mode $(0, 5)$ (Panels (a)--(b)) and mode $(-2, 3)$ (Panels (c)--(d)) of the quasi-geostrophic (QG) ocean model. The black and the solid blue curves are the truth and the posterior mean estimate. The shading area indicates one standard deviation (S.D.) from the posterior mean. The light blue curves are the sampled realizations of the corresponding mode from the posterior distribution. }
\end{figure}

Figure \ref{fig:smoother} shows the result of the state estimation (the posterior distribution) of the ocean field using the Lagrangian data assimilation. Data assimilation is applied to the spectral form of the model and therefore the results in Figure \ref{fig:smoother} are the posterior distribution of two specific spectral modes. The posterior estimates of other spectral modes have qualitatively similar behavior. The solid blue curve is the posterior mean estimate, which is usually used as the optimal point estimate for deterministic eddy identification methods. However, it is notable that the posterior mean can be different from the truth (black curve), which can be regarded as one random realization from the underlying turbulent model. The shading area indicates one standard deviation from the posterior mean. The difference in various sampled realizations, marked by light blue curves, is due to such uncertainty in the state estimation.

\subsection{\label{sec:owuncertainty}The effect of uncertainty on eddy diagnostics}

\begin{figure}
    \begin{center}
        \includegraphics[]{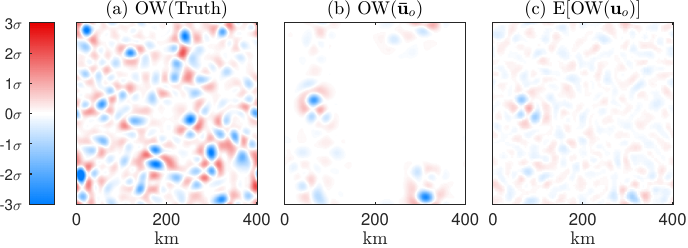}
    \end{center}
    \caption{\label{fig:evofow} Comparison between applying the OW parameter using the true flow field (Panel (a)), using the posterior mean ocean state, $\mathrm{OW}(\bar{\mathbf{u}}_o)$, (Panel (b)), and the expected value of the OW parameter based on the multiple realizations of the ocean field sampled from the posterior distribution, $\mathbb{E}[\mathrm{OW}(\mathbf{u}_o)]$, (Panel (c)). The result in Panel (b) can be regarded as one of the standard deterministic eddy identification diagnostics. The result in Panel (c) considers the uncertainty in the state estimation of the ocean. The same contour level sets are used on the two panels. To highlight the role of strong uncertainty in affecting the eddy identification skill, the distribution of the estimated ocean state in this experiment is obtained from the Lagrangian data assimilation using only four observed ice floe trajectories. This is different from using 40 floes as in all other experiments.
    }
\end{figure}

Figure \ref{fig:evofow} shows the difference between applying the OW parameter to the posterior mean ocean state, $\mathrm{OW}(\bar{\mathbf{u}}_o)$, (Panel (b)) and the expected value of the OW parameter based on the multiple realizations of the ocean field sampled from the posterior distribution, $\mathbb{E}[\mathrm{OW}(\mathbf{u}_o)]$, (Panel (c)). The former can be regarded as one of the standard deterministic eddy identification diagnostics based on the optimal point estimate (e.g., the posterior mean). The latter considers the uncertainty in the state estimation of the ocean. To highlight the role of strong uncertainty in affecting the eddy identification skill, the distribution of the estimated ocean state in this experiment is obtained from the Lagrangian data assimilation using only four observed ice floe trajectories. The results from these two methods are significantly contrasting to each other. A much fewer number of eddies are identified in Panel (b). This is because when the uncertainty dominates the state estimation due to the lack of sufficient observations, the posterior mean estimate of the ocean field will behave similarly to the model equilibrium mean value. Since the model is an anomaly model, the mean value is zero. Therefore, the posterior mean severely underestimates the amplitude of the ocean field and affects the subsequent eddy identification skill. Note that when 40 floe observations are utilized, the OW parameter based on the posterior mean ocean state remains to have a significant error by identifying only 50\% of the eddies (see the following subsections).

The effect of uncertainty on computing the OW parameter can be demonstrated theoretically. Since the estimated ocean state is characterized in a probabilistic way by a PDF, it can be regarded as a random variable. Therefore, a mean-fluctuation decomposition can be applied to the estimated ocean state $\mathbf{u}_o = \bar{\mathbf{u}}_o + \mathbf{u}_o'$, where $\bar{\mathbf{u}}$ is the deterministic mean component representing the optimal point estimate and $\mathbf{u}_o'$ is the random fluctuation component. Denote by $\mathrm{OW}(\bar{\mathbf{u}}_o)$ the OW parameter based on the posterior mean ocean state (Panel (a) in Figure \ref{fig:evofow}) and $\mathbb{E}[\mathrm{OW}(\mathbf{u}_o)]$ the expected value of the OW parameter based on the multiple realizations of the ocean field sampled from the posterior distribution (Panel (b) in Figure \ref{fig:evofow}). In light of the definition of the OW parameter \eqref{eqn:OW1}--\eqref{eqn:OW2}, the relationship between $\mathrm{OW}(\bar{\mathbf{u}}_o)$ and $\mathbb{E}[\mathrm{OW}(\mathbf{u}_o)]$ is given by
\begin{equation}
    \label{eqn:evofow}\mathbb{E}[\mathrm{OW}(\mathbf{u}_o)] = \mathrm{OW}(\bar{\mathbf{u}}_o) + \underbrace{\mathbb{E}[(u_x')^2] - 2\mathbb{E}[u_x'v_y'] + \mathbb{E}[(v_y')^2] + 4\mathbb{E}[v_x'u_y']}_{\text{Additional terms due to uncertainty.}},
\end{equation}
where $\cdot_x$ and $\cdot_y$ represent the spatial derivative with respect to $x$ and $y$ respectively. For notation simplicity, the two components of $\mathbf{u}_o$, namely $(u_o,v_o)$, are written as $(u,v)$ in \eqref{eqn:evofow} and in \ref{sec:evofow} that includes the detailed derivations. The additional terms on the right-hand side of \eqref{eqn:evofow} involve the contribution from the fluctuation part, which results from sampling the ocean field from the PDF of the estimated state. In general, as most of the eddy identification criteria are nonlinear with respect to the flow field \cite{branicki2011lagrangian, mancho2013lagrangian, vortmeyer2016detecting, vortmeyer2019comparing, d2009comparison, souza2011comparison}, changing the order of taking the expectation and applying the eddy diagnostic will have a significant impact on the diagnostic results. In addition, the fluctuation part will have a significant contribution to the eddy diagnostic when the expectation is taking over the identified eddy field based on the sampled realizations of the ocean from its probabilistic state estimation.

\subsection{Probabilistic quantification of the number of eddies}
Due to the uncertainty in the state estimation of the ocean field, it is natural to provide an uncertainty quantification when counting the number of eddies.

\begin{figure}
    \begin{center}
        \includegraphics{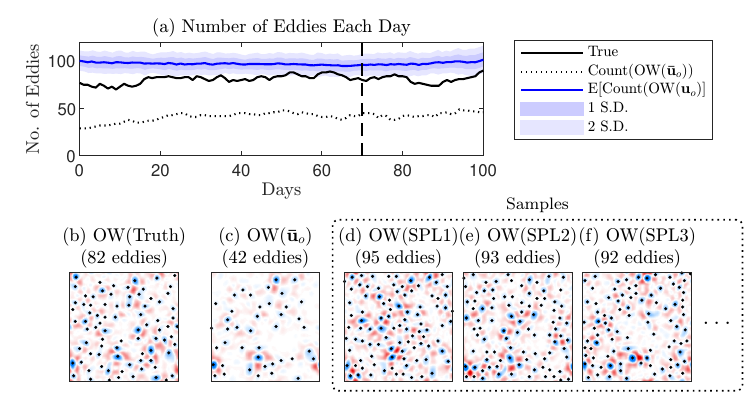}
    \end{center}
    \caption{\label{fig:count} The number of identified eddies. Panel (a): the number of identified eddies in the true ocean field as a function of time (black solid curve), the number of eddies identified using the optimal point estimate, i.e., the posterior mean, of the ocean $\mbox{OW}(\bar{\mathbf{u}}_o)$ (black dashed curve), and the number of eddies identified using the probabilistic method $\mbox{OW}({\mathbf{u}}_o)$. The probability of the number of identified eddies in each day is presented by the expectation value, $\mathbb{E}[\mbox{Count}((\mbox{OW}({\mathbf{u}}_o))]$ (blue solid curve), and the standard deviations (S.D.) (blue shading areas). Panel (b): the identified eddies from the true flow field on day 70. Panel (c): the identified eddies using the deterministic method, based on the posterior mean field, namely $\mbox{OW}(\bar{\mathbf{u}}_o)$, on day 70. Panels (d)--(f): the identified eddies from three sampled ocean realizations (named SPL 1-3) drawn from the posterior distribution at day 70. The black dots in Panels (b)--(d) represent the eddy cores.}
\end{figure}

Figure \ref{fig:count} illustrates the number of the identified eddies by applying the probabilistic eddy identification method described in Section \ref{sec:procedure} and Figure \ref{fig:overview}, which is compared with the deterministic method based on the optimal point estimate, namely the posterior mean field. It is seen that using the posterior mean as a proxy for the true ocean flow severely underestimates the number of eddies compared with the truth. Only about 50\% of the eddies are identified. This is again due to the underestimation of the flow amplitude by omitting the information in the fluctuation part. In contrast, the number of eddies identified using the sampled ocean realizations drawn from the posterior distribution of the estimated ocean field is much closer to the truth. Note that there is still a bias in the eddy count (the gap between the black and blue curves). Such a bias is, in fact, introduced by the approximations in the approximate linear stochastic model used in the data assimilation scheme. See \ref{sec:twinexperiment} for a perfect model twin experiment, which explains such a model error, which often appears in practice, in affecting the eddy identification skill. Nevertheless, even in the presence of model error in applying data assimilation for estimating the ocean state, the probabilistic eddy count outperforms the deterministic point estimate method. Note that as the number of observations increases, the bias in the probabilistic method reduces, and the number of the identified eddies eventually converges to the truth.

\subsection{Probabilistic occurrence of each eddy}
The probabilistic framework can also provide a likelihood of the occurrence of each eddy. This is achieved by selecting a date and a location and then checking the OW parameter associated with each sampled ocean realization. If an eddy core is identified within a certain radius of the prescribed location (in this study, 10 km), then that eddy is identified in the sample. Here, an eddy core is defined as a local minimum of the OW parameter, which is below the cutoff threshold. The ratio of the total number counted over the number of sampled realizations leads to the probability of the occurrence of such an eddy.

\begin{figure}
    \begin{center}
        \includegraphics{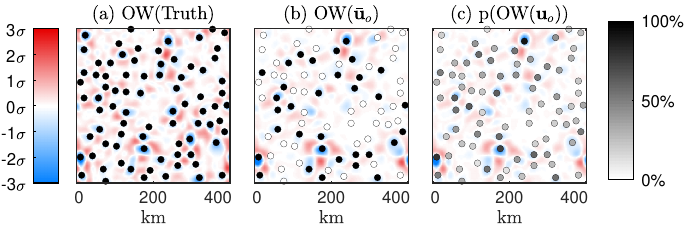}
    \end{center}
    \caption{\label{fig:identification} Probabilistic occurrence of each eddy.
    Panel (a): identified eddies based on the true ocean field at day 70, where black dots mark the identified eddy cores. Panel (b): identified eddies based on the optimal point estimate, i.e., the posterior mean, of the ocean $\mbox{OW}(\bar{\mathbf{u}}_o)$. Black dots mark the eddies that appear in the truth and are identified. The eddies that are not identified are marked by hollows. Panel (c): the probability of each eddy identified using the probabilistic method. The color of the dots shows the probability.
    }
\end{figure}

Panel (a) of Figure \ref{fig:identification} shows identified eddies based on the true ocean field at day 70, where black dots mark the identified eddy cores. Panel (b) illustrates the identified eddies based on the optimal point estimate, i.e., the posterior mean, of the ocean $\mbox{OW}(\bar{\mathbf{u}}_o)$. Only about half of the eddies are identified using such a deterministic method. This is not surprising as the deterministic criterion naturally leads to large biases when the estimated ocean state contains considerable uncertainty. In contrast, Panel (c) shows the results using the probabilistic method. The identification results are characterized by probabilities, which indicate the likelihood of the occurrence of each eddy. It is worth highlighting that, to provide an intuitive outcome in practice, a deterministic conclusion can nevertheless be derived from the probabilistic eddy identification results. To this end, a threshold (e.g., 50\%) can be applied to the probabilities and only keep the identified eddies above such a threshold representing a certain confidence level. Note that the result from such a probabilistic-induced deterministic conclusion will still be fundamentally different from the one shown in Panel (b) since the nonlinearity in the eddy identification criteria considers the contribution from the fluctuations that provide additional information to help detect eddies using the probabilistic method.

\subsection{\label{sec:lifetime}Probabilistic quantification of eddy lifetime and size}

In addition to counting the number of eddies or calculating the probability of the eddy occurrence at a fixed time instant, there exist many other important quantities that are related to eddy dynamics. Among these quantities, eddy lifetime and eddy size are the two widely used indicators that characterize the underlying physics of the eddies. In this study, the eddy lifetime is calculated from a time series of ocean states by tracking the eddy core over time. To calculate the lifetime of a specific eddy, an initial specification of the eddy location and at a given date is used. If an eddy core is identified at this date within a certain initial search distance of the specified location, then the eddy lifetime can be calculated. Between consecutive snapshots of the OW parameter, eddy cores can be matched between the snapshots if they are within a certain tracking distance of each other. The eddy core is tracked forward and backward in time until it no longer meets the eddy criterion. This can be because the OW parameter passes above the cutoff threshold or because the core is no longer a local minimum. The duration of time that the eddy exists within the snapshots is used to calculate the eddy lifetime. In this study, day 70 is used as the initial specification of the eddy, and a distance of 10 km between two consecutive days is used. On the other hand, eddy size is calculated from a single snapshot of the ocean state. According to the criterion used in this study, the eddy boundary is defined as the closed contour of the OW parameter at the cutoff threshold containing the eddy core. If an eddy is identified according to the eddy criterion, the area enclosed by the eddy boundary can be calculated. Notably, when quantities such as the eddy lifetime are computed, sampling the ocean realizations at a fixed time instant is insufficient as eddies evolve in time. Therefore, the sampling is taken from the joint posterior distribution of the state estimation at different time instants. Details of such a sampling algorithm are included in \ref{sec:cg}.

\begin{figure}
    \begin{center}
        \includegraphics[width=\textwidth]{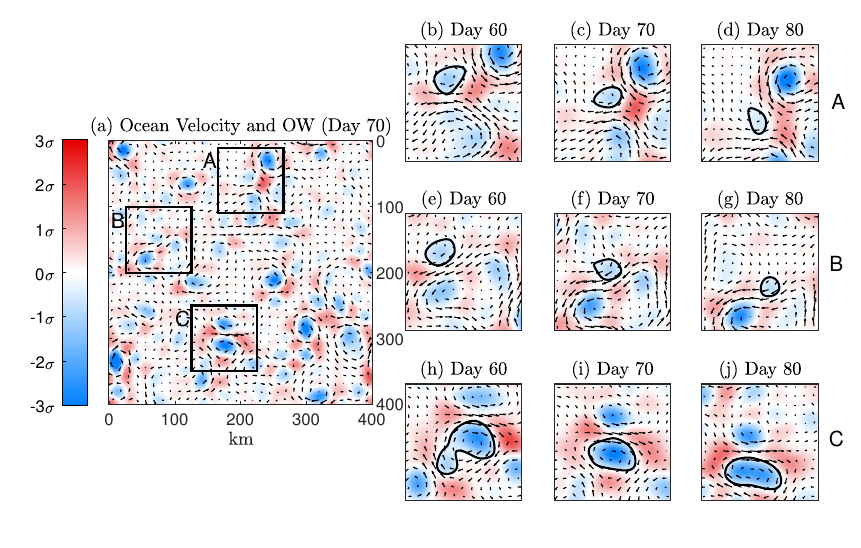}
    \end{center}
    \caption{\label{fig:eddysnapshot} Eddies for case studies. Panel (a): the true ocean field at day 70 and the three eddies, labeled ``A'', ``B'', and ``C'', used as case studies in Figures \ref{fig:combined} and \ref{fig:lifetime_snapshots}. For illustration purpose, the OW parameter has been normalized to have a standard deviation of 1. Panels (b)-(d), (e)-(g), and (h)-(j): snapshots of the enlarged square areas containing eddies A, B, and C, marked by closed cycles in black, at different times.}
\end{figure}

\begin{figure}
    \begin{center}
        \includegraphics{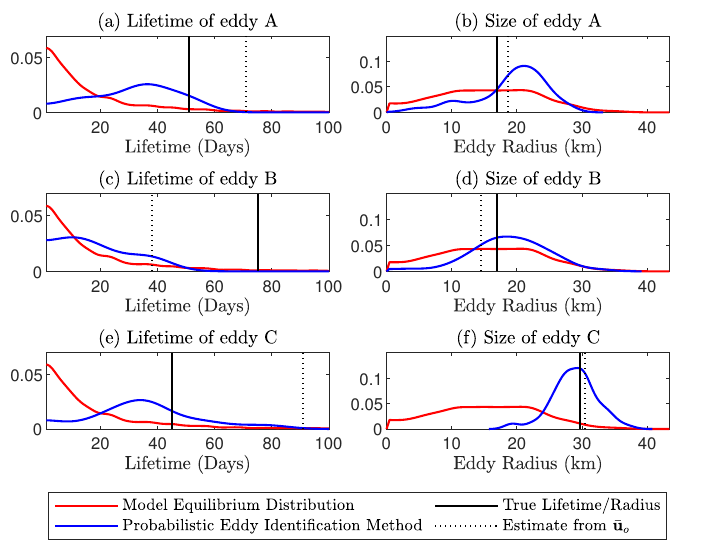}
    \end{center}
    \caption{\label{fig:combined} Distributions of eddy lifetime and eddy size for the three eddies highlighted in figure \ref{fig:eddysnapshot}. Panels (a), (c), and (e): the distributions of eddy lifetime for eddies A, B, and C. The red curve represents the distribution of eddy lifetimes by collecting all the eddies simulated from the true two-layer QG model. Such a distribution reflecting the climatological properties of the model that is in the absence of observations. The blue curve shows the distribution of eddy lifetime for the corresponding eddy using the probabilistic eddy identification framework based on an ensemble of sampled ocean realizations. The solid black line shows the true eddy lifetime for the corresponding eddy. The dotted line shows the lifetime of the eddy computed based on the recovered ocean field using the posterior mean estimate. Panels (b), (d), and (f): the distributions of size for the three eddies.}
\end{figure}

Figure \ref{fig:eddysnapshot} illustrates the time evolution of three eddies, which are utilized as case studies to compare the identification of the eddy lifetime and size using the probabilistic framework and its deterministic counterpart. The red curve in the left column of Figure \ref{fig:combined} represents the overall distribution of eddy lifetimes by collecting all the eddies simulated from a long simulation using the true two-layer QG model. It reflects the climatological properties of the model in the absence of observations. The distribution indicates that the lifetime of eddies in the flow field varies significantly. While many eddies are transient, with a lifetime that is shorter than 20 days, others can last much longer. The eddies selected as case studies belong to the latter category. The solid black lines in the left column of Figure \ref{fig:combined} show the true lifetime of the three eddies marked in Figure \ref{fig:eddysnapshot}. Eddy A and Eddy C last for 45 to 50 days, which is longer than most of the eddies in the field, while Eddy B has an extremely long lifetime of about 75 days. The dashed black curve in each panel shows the estimated lifetimes of these eddies computed from the estimated ocean field using the posterior mean estimate from the Lagrangian data assimilation. Estimating the eddy lifetime using a deterministic method leads to a large error. For Eddy A, the estimated lifetime is longer than 70 days, while the truth is 50 days. For Eddy B and Eddy C, the errors in the estimation are even more severe. The estimated lifetime of the long-lasting Eddy B is only half of the true value, while the lifetime estimation of Eddy C is nearly doubled to 90 days. The significant biases in the estimation primarily come from the large discrepancy between the truth and the posterior mean estimate of the ocean, which is the optimal point estimate, due to ignorance of the considerable uncertainty in state estimation. This again highlights the necessity of incorporating uncertainty in eddy identification, even for computing the mean estimated lifetime, since the eddy diagnostic is a nonlinear function of the estimated states. In addition to the absolute difference between the estimation and the true lifetime, the lack of quantification of estimation accuracy is another major shortcoming, especially since the eddy lifetime is affected by the turbulent nature of the underlying flow field.

In contrast, the probabilistic eddy identification framework characterizes the eddy lifetime fundamentally differently by providing a distribution for the lifetime estimation. First, this distribution differs significantly from the one based on the model climatology (red curve). This indicates that observational information plays a significant role in reducing the uncertainty in estimating the ocean field and the eddy lifetime. Second, the distribution provides a better likelihood estimation of the lifetime for Eddy A and Eddy C. The true value of the eddy lifetime is not only close to the peak of the distribution but also within the range of the quantified uncertainty in the estimation. Although the probabilistic framework still gives only a small probability for Eddy B, the entire distribution from the framework (blue curve) shifted towards the larger value compared with the model climatology. This means the probabilistic framework implies a longer lifetime of such an event than most other events, indicating the possibility of this as an extreme event. Similarly, the eddy size distributions computed from the probabilistic eddy identification framework are shown in the right column of Figure \ref{fig:combined}, illustrating a significant reduction of the uncertainty related to the model climatology. The truth is also located near the peak of the distribution. Note that the estimated value using the deterministic method (dashed black line) is closer to the truth for inferring the eddy size than the eddy lifetime. This is because the eddy size is computed from the recovered ocean field at a single time instant while the lifetime is estimated from a sequence of snapshots, where uncertainty is propagated and enlarged in the latter case.

\begin{figure}
    \begin{center}
        \includegraphics{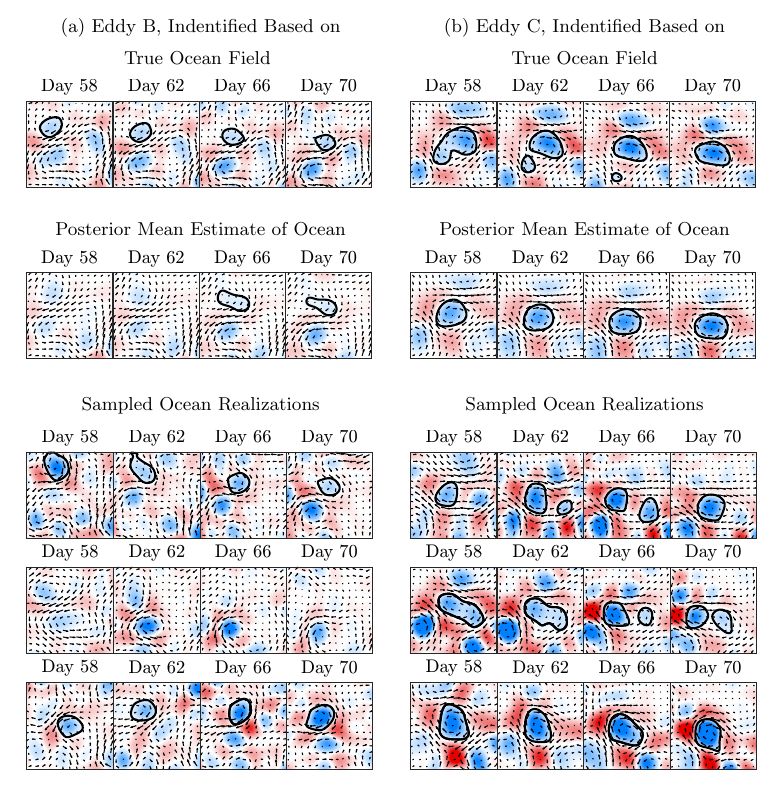}
    \end{center}
    \caption{\label{fig:lifetime_snapshots}Snapshots of the time evolution of eddies B and C identified based on the true ocean flow, the estimated ocean flow using the posterior mean estimate, and three sampled trajectories from the posterior distribution of the ocean state. Panels (a): the time evolution of Eddy B. Panels (b): the time evolution of Eddy C. }
\end{figure}

To better understand the reasoning of the difference between the probabilistic and deterministic methods in Figure \ref{fig:combined}, Figure \ref{fig:lifetime_snapshots} includes snapshots of the time evolution of Eddy B and Eddy C identified based on the true ocean flow, the estimated ocean flow using the posterior mean estimate, and three sampled trajectories from the posterior distribution of the ocean state. The second row in Panel (a) shows that the magnitude of the OW parameter is underestimated based on the recovered ocean field using the posterior mean estimate. Therefore, Eddy B is not identified on days 58 and 62, which leads to a severe underestimation of the eddy lifetime. In contrast, although the eddy is not detected in some of the sampled realizations of the ocean, there is a high chance the eddies are identified through the period, increasing the likelihood of the eddy with a longer lifetime from the probabilistic eddy identification framework. On the other hand, in the first row of Panel (b), due to the highly turbulent ocean field in the local area, an eddy that appears on day 58 splits into two eddies. However, since the posterior mean estimate of the ocean field smooths out the fluctuations, the associated eddies are more regular, as is seen in the second row. In contrast, when the uncertainty in the state estimation of the ocean is taken into account, the eddies associated with the sampled ocean realizations exhibit similar complex dynamical features as the truth, where eddies are more prone to split and have irregular shapes, as is seen in the last three rows in Panel (b).

\section{\label{sec:conclusions} Conclusions}

Due to the intrinsic strong nonlinearity in the eddy identification diagnostics, even a small uncertainty in estimating the flow field can cause a significant error in the identified eddies. In this paper, a general probabilistic eddy identification framework is developed to highlight the importance of considering the uncertainty propagation through the nonlinear diagnostics for eddy identification problems. The framework can adapt existing identification methods to incorporate uncertainty into the diagnostic. The probabilistic eddy identification framework is based on sampling an ensemble of ocean realizations from the probabilistic state estimation. The corresponding eddy statistics are aggregated from the diagnostic results by applying traditional eddy diagnostics to individual realizations. The framework is applied to a scenario mimicking the Beaufort Gyre marginal ice zone, where large uncertainty appears in estimating the ocean field using Lagrangian data assimilation with sparse ice floe trajectories.  The skills in counting the number of eddies and computing the probability of each eddy event are significantly improved under the probabilistic framework. Notably, accounting for the nonlinear propagation of uncertainty in diagnostics provides a more accurate mean estimate than standard deterministic methods in estimating eddy lifetime. It also facilitates uncertainty quantification in inferring such a crucial dynamical quantity.

Several research topics can be further addressed based on the probabilistic framework developed in this paper. First, the OW parameter is utilized as the eddy diagnostic criterion in this paper due to its simple form. However, it has been shown that such an Eulerian approach may not always provide the most skillful solutions. There are many Lagrangian eddy identification methods \cite{branicki2011lagrangian, nolan2020finite, mancho2013lagrangian, vortmeyer2016detecting, vortmeyer2019comparing}, which incorporate dynamical information from the ocean and are potentially more robust. Therefore, one natural future direction is to compare these two types of approaches under the probabilistic framework and explore how the uncertainty in state estimation affects the results. Second, exploring the impact of each source of uncertainty is a crucial topic. It helps to optimally determine the effort to put into improving the measurements and the models that advance eddy identification. Third, some simplifications were made when applying the probabilistic framework to the scenario mimicking the Beaufort Gyre marginal ice zone. These simplifications are adopted to facilitate the simulations and allow the focus on comparing the probabilistic framework with its deterministic counterpart. They can easily be removed when studying the real-world scenario. As the Arctic Ocean becomes an important indicator of climate change, applying the probabilistic framework to the newly developed real observational data set \cite{watkins2024observing, lopez2021library} to reveal the variation of the eddy statistics over the past few decades has practical significance.

\appendix

\section{Details of the sea ice floe model \label{sec:ice}}

Ice floes are formed when ice sheets break up. Especially during the summer, ice floes can be free-floating, driven by the local characteristics of the ocean. To utilize satellite observations of ice floe trajectories for data assimilation, it is necessary to incorporate a model that can account for intricate ice floe dynamics. In this section, we present a DEM model of ice floe dynamics that captures complex sea ice and ice-ocean interactions, which can be used in data assimilation applications despite its high dimensionality and nonlinearity.

The ice floes are subject to ocean drag forces in a one-way interaction where the ocean drag force on the ice floes is calculated using a quadratic drag approximation. Further, the ocean forces and torques acting on the floes are assumed to be uniform over the shape of the floe, allowing for the forces to be explicitly calculated without the need to calculate surface integrals in the floe dynamics. Finally, the contact forces between floes can be approximated by using cylindrical floe shapes. Despite these simplifications, this model can capture many of the rich features of sea ice dynamics necessary to account for in a real observational setting.

Consider an ice floe model with arbitrary shape. The governing equation of motion is given by:
\begin{equation} \label{eqn:floe}
\begin{split}
    \frac{\d\mathbf{x}^\ell}{\d t} =& \mathbf{u}^\ell\\
    m_\ell \frac{\d\mathbf{u}^\ell}{\d t} =& \int_{A_\ell} \mathbf{F}_\mathrm{o}^\ell(\mathbf{x}) \, d\mathbf{x}  \\
    \frac{\d\Omega^\ell}{\d t} =& \omega^\ell  \\
    I_\ell \frac{\d\omega^\ell}{\d t} =& \int_{A_\ell} \tau_\mathrm{o}^\ell(\mathbf{x}) \|\mathbf{x} - \mathbf{x}^\ell\|^2 \, d\mathbf{x},
\end{split}
\end{equation}
where the superscript $\ell$ represents the $\ell$-th floe. In \eqref{eqn:floe}, $\mathbf{x}^\ell$ and $\Omega^\ell$ are the displacement and the angular displacement of the floe while $\mathbf{u}^\ell$ and $\omega^\ell$ are the velocity and angular velocity, respectively.
The forcing and the torque are given by
\begin{equation}
    \mathbf{F}_\mathrm{o}^\ell (\mathbf{x}, t) = d_\mathrm{o} \rho_\mathrm{o} \left(\mathbf{u}_\mathrm{o}(\mathbf{x}, t) - \mathbf{u}^\ell\right)\left|\mathbf{u}_\mathrm{o}(\mathbf{x}, t) - \mathbf{u}^\ell\right|
\end{equation}
and
\begin{equation}
    \tau_\mathrm{o}^\ell(\mathbf{x}, t) = d_\mathrm{o} \rho_\mathrm{o} \left(\frac{\nabla \times \mathbf{u}_\mathrm{o}(\mathbf{x}, t)}{2} - \omega\right)\left|\frac{\nabla \times \mathbf{u}_\mathrm{o}(\mathbf{x}, t)}{2} - \omega\right|,
\end{equation}
respectively.
Define the floe area and the second polar moment of area as
\begin{equation}
    A = \int_A \, \d\mathbf{x} \qquad\mbox{and}\qquad J = \int_A r^2 \, \d\mathbf{x}.
\end{equation}
Then the mass $m_\ell$ and the moment of inertia $I_\ell$ are given by
\begin{equation}
    m_\ell = \int_{A_\ell} \rho_\mathrm{ice} h_\ell \, \d\mathbf{x} = \rho_\mathrm{ice} A h_\ell  \qquad\mbox{and}\qquad I_\ell = \int_{A_\ell} h_\ell \rho_\mathrm{ice} r^2 \, \d\mathbf{x} = \rho_\mathrm{ice} J h_\ell.
\end{equation}

\section{Details of the quasi-geostrophic ocean model\label{sec:qg}}

The ocean model is a two-layer Quasi-Geostrophic (QG) model with periodic boundary conditions on a square domain. The ocean state is characterized by the stream functions $\psi_i(x, y)$ and potential vorticities (PV) $q_i(x, y)$ of each layer $i = 1, 2$. The level curves of the stream function, $\psi_i$, correspond to streamlines of the velocity field, which guarantees an incompressible flow. The ocean velocity field for each layer can thus be calculated as
\begin{equation}
(u_i, v_i) = \left(-\frac{\partial \psi_i}{\partial y}, \frac{\partial \psi_i}{\partial x}\right), \quad i = 1, 2. \label{eqn:vel}
\end{equation}

The formulation of the QG equations follows the version in \cite{arbic2004baroclinically}. The PDEs which govern the time evolution of $\psi_i$ and $q_i$ are as follows:
\begin{align}
\frac{\partial q_1}{\partial t} + \overline{u_1} \frac{\partial q_1}{\partial x}
+ \frac{\partial \overline{q}_1}{\partial y}\frac{\partial \psi_1}{\partial x} + J(\psi_1, q_1) =& \mathrm{ssd} \\
\frac{\partial q_2}{\partial t} + \overline{u_2} \frac{\partial q_2}{\partial x}
+ \frac{\partial \overline{q_2}}{\partial y} \frac{\partial \psi_2}{\partial x} + J(\psi_2, q_2) =& -R_2 \nabla^2 \psi_2 + \mathrm{ssd}.
\end{align}
Here ``ssd'' represents small-scale dissipation, which are higher-order derivative terms that are ignored. The term $J$ is the Jacobian
\begin{equation}
J(\psi, q) = \frac{\partial \psi}{\partial x} \frac{\partial q}{\partial y} - \frac{\partial \psi}{\partial y} \frac{\partial q}{\partial x}.
\end{equation}
The stream functions further satisfy
\begin{align}
q_1 =& \nabla^2 \psi_1 + \frac{(\psi_2 - \psi_1)}{(1 + \delta) L_d^2} & q_2 =& \nabla^2 \psi_2 + \frac{\delta(\psi_1 - \psi_2)}{(1 + \delta) L_d^2}.
\end{align}
where $\delta = H_1/H_2$, $H_i$ is the depth of each layer, and $L_d$ is the deformation radius.
The terms $\partial \overline{q}_1/\partial y$ and $\partial \overline{q}_2/\partial y$, despite the notation, are parameters representing the mean PV gradients for each layer and are given by
\begin{align}
\frac{\partial \overline{q_1}}{\partial y} =& \frac{\overline{u_1} - \overline{u_2}}{(1 + \delta)L_d^2} & \frac{\partial \overline{q_2}}{\partial y} =& \frac{\delta (\overline{u_2} - \overline{u_1})}{(1 + \delta) L_d^2}
\end{align}
where $\overline{u_1}$ and $\overline{u_2}$ are the mean ocean velocities. The final parameter, $R_2$, is the  decay  rate  of  the  barotropic  mode
\begin{equation}
R_2 = \frac{f_0 d_{\mathrm{Ekman}}}{2H_2}
\end{equation}
where $f_0$ is the Coriolis parameter and $d_{\mathrm{Ekman}}$ is the bottom boundary layer thickness. Note that in this formulation we use a constant Coriolis force throughout the domain.

The parameters in the floe and the ocean models that generate the true signal are summarized in Table \ref{tab:parameters}.

\begin{table}
    \centering
    \begin{tabular}{cc}
        Parameter & Value \\ \hline
        Ocean density & $\rho_{\text{ocn}} =1000$kg/m$^3$  \\
        Ice density &$\rho_{\text{ice}}=920$kg/m$^3$ \\
        Ocean drag coefficient & $c_{\text{ocn}}=3\times 10^{-3}$\\
        Top layer mean ocean velocity & $\overline{u_1}=2.58$km/day\\
        Bottom layer mean ocean velocity & $\overline{u_2}=1.032$km/day\\
        Top layer mean potential vorticity&$\frac{\partial \overline{q_1}}{\partial y}=0.0265$km$^{-1}$day$^{-1}$\\
        Bottom layer mean potential vorticity&$\frac{\partial \overline{q_2}}{\partial x}=-0.0212$km$^{-1}$day$^{-1}$\\
        Coupling parameter &$R_1=6.9\times10^{-5}$km$^{-1}$\\
        Decay rate of the barotropic mode & $R_2=1$day$^{-1}$\\
        Deformation radius	&$L_d=5.7$km\\
        Ratio of upper-to lower-layer depth & $\delta=0.8$\\
        OW parameter cutoff & $0.2 \sigma_{\mathrm{OW}}$ \\
        Domain size & 400 km $\times$ 400 km \\
        Number of ice floes & 40 floes \\
        Floe thickness & 2 m \\
        Trajectory length & 100 days \\
        Observational noise in location & 250 m \\
    \end{tabular}
    \caption{Parameters in the ice-ocean model.}
    \label{tab:parameters}
\end{table}

\section{Details of the stochastic ocean model \label{sec:ocean}}

To facilitate data assimilation, a stochastic forecast model of the ocean is adopted to approximate the forecast statistics from the two-layer QG model. Although the stochastic model cannot capture every trajectory of the original QG model due to its turbulent nature, the forecast statistics can be approximated quite accurately. Such a stochastic forecast model has been widely used to accelerate data assimilation \cite{majda2016introduction, farrell1993stochastic, berner2017stochastic, harlim2008filtering, branicki2013non}.

The stochastic forecast model is building upon each Fourier model. The stochastic model describes incompressible flows, as the QG system, and is periodic on the domain $[0, x_\mathrm{max})^2$. Denote by $\hat{u}_\mathbf{k}(t)$ the $\mathbf{k}$-th Fourier coefficient, where $\mathbf{k} = (k_1, k_2)^\mathsf{T} \in \mathbb{Z}^2$ is the two-dimensional wave number and the coefficients satisfy $\hat{u}_{-\mathbf{k}} = \overline{\hat{u}_\mathbf{k}}$, the complex conjugate.
The real-valued stream function, $\psi$, is given by
\begin{equation}
    \psi(\mathbf{x}, t) = \sum_{\mathbf{k} \in \mathbb{Z}^2} \hat{u}_\mathbf{k}(t) \exp\left(-\frac{2\pi i}{x_\mathrm{max}}(\mathbf{k}\cdot\mathbf{x})\right)
\end{equation}
where the level sets of the stream function are the streamlines of the velocity field, so the corresponding ocean velocity is given by $\mathbf{u}_\mathrm{o} = (u_\mathrm{o}, v_\mathrm{o})^\mathsf{T} = \left(\frac{\partial \psi}{\partial y}, -\frac{\partial \psi}{\partial x}\right)$.
From the Fourier coefficients, the ocean velocity is reconstructed as
\begin{equation}
    \label{eqn:ocean_velocity}
    \mathbf{u}_\mathrm{o}(\mathbf{x}, t) = \sum_{\mathbf{k} \in \mathbb{Z}^2} \hat{u}_\mathbf{k}(t) \exp\left(-\frac{2\pi i}{x_\mathrm{max}}(\mathbf{k}\cdot\mathbf{x})\right) \mathbf{r}_\mathbf{k}
\end{equation}
with $\mathbf{r}_\mathbf{k} = \left(\frac{2\pi ik_2}{x_\mathrm{max}}, -\frac{2\pi ik_1}{x_\mathrm{max}}\right)^\mathsf{T}$ is the eigenvector that represents the relationship between the two velocity component. The incompressibility is reflected in the eigenvector.
In addition, the ocean vorticity, $\nabla \times \mathbf{u}_\mathrm{o} = \frac{\partial v_\mathrm{o}}{\partial x} - \frac{\partial u_\mathrm{o}}{\partial y}$, is given by
\begin{equation}
    \label{eqn:ocean_vorticity}
    \nabla \times \mathbf{u}_\mathrm{o}(\mathbf{x}, t) = \sum_{\mathbf{k} \in \mathbb{Z}^2} \hat{u}_\mathbf{k}(t) \exp\left(-\frac{2\pi i}{x_\mathrm{max}}(\mathbf{k}\cdot\mathbf{x})\right) \left(\frac{2\pi |\mathbf{k}|}{x_\mathrm{max}}\right)^2.
\end{equation}

Each Fourier coefficient $\hat{u}_\mathbf{k}(t)$ is a function of time $t$. Here the stochastic ocean model is used to describe the Fourier coefficients. Each Fourier coefficient is governed by an independent Ornstein-Uhlenbeck (OU) process \cite{gardiner2004handbook}
\begin{equation}
    \label{eqn:ouprocess} \frac{\d\hat{u}_\mathbf{k}}{\d t} = (-a_\mathbf{k} + \phi_\mathbf{k} i) \hat{u}_\mathbf{k} + f_\mathbf{k} + \sigma_\mathbf{k} \dot{W}_\mathbf{k}(t)
\end{equation}
with $\dot{W}_\mathbf{k}$ is a complex-valued Gaussian white noise and the parameters $\sigma_\mathbf{k} > 0$, $a_\mathbf{k}$ and $\phi_\mathbf{k}$ are all real-valued while $f_\mathbf{k}$ is complex.
Note that even though the OU processes for the Fourier coefficients given in equation \eqref{eqn:ouprocess} are independent, the resulting velocity field in physical space still has spatial correlations.
In order to ensure that the ocean velocity is real-valued it is necessary that $\hat{u}_{-\mathbf{k}} = \overline{\hat{u}_\mathbf{k}}$, so the parameters and noise have the additional restriction that $a_{-\mathbf{k}} = a_\mathbf{k}$, $\phi_{-\mathbf{k}} = -\phi_\mathbf{k}$, $f_{-\mathbf{k}} = \overline{f_\mathbf{k}}$, and $\sigma_{-\mathbf{k}} = \sigma_\mathbf{k}$ with $\dot{W}_{-\mathbf{k}} = \overline{\dot{W}_\mathbf{k}}$.

The equilibrium mean, covariance, and decorrelation time of each independent OU process given in \eqref{eqn:ouprocess} can be written explicitly in terms of the parameters $a_\mathbf{k}$, $\phi_\mathbf{k}$, $f_\mathbf{k}$, and  $\sigma_\mathbf{k}$ as \cite{chen2023uncertainty}
\begin{equation}
    \label{eqn:eqstats}
    \operatorname{Mean}(\hat{u}_\mathbf{k}) = \frac{f_\mathbf{k}}{a_\mathbf{k} - i\phi_\mathbf{k}}, \qquad
    \operatorname{Var}(\hat{u}_\mathbf{k}) = \frac{\sigma_\mathbf{k}^2}{2a_\mathbf{k}}, \qquad\mbox{and}\qquad
    T_\mathrm{corr} = \frac{1}{a_\mathbf{k} - i\phi_\mathbf{k}},
\end{equation}
where the decorrelation time is defined as
\begin{equation}
    T_\mathrm{corr} = \int_0^\infty \frac{\mathbb{E}[(\hat{u}_\mathbf{k}(0) - \operatorname{Mean}(\hat{u}_\mathbf{k}))(\hat{u}_\mathbf{k}(t) - \operatorname{Mean}(\hat{u}_\mathbf{k}))]}{\operatorname{Var}(\hat{u}_\mathbf{k})} \, \d t.
\end{equation}
Inversely, the equations given in \eqref{eqn:eqstats} can be inverted to write the parameters in terms of a given a set of equilibrium statistics
\begin{align}
    a_\mathbf{k} =& \operatorname{Re}
\left[\frac{1}{T_\mathrm{corr}}\right] &  \phi_\mathbf{k} =& -\operatorname{Im}
\left[\frac{1}{T_\mathrm{corr}}\right] \\ f_\mathbf{k} =& \frac{\operatorname{Mean}(\hat{u}_\mathbf{k})}{T_\mathrm{corr}} & \sigma_\mathbf{k} =& \sqrt{2 \operatorname{Var}(\hat{u}_\mathbf{k}) \operatorname{Re}
\left[\frac{1}{T_\mathrm{corr}}\right]}.
\end{align}
In this way, the parameters can be systematically calibrated to a specified equilibrium mean, variance, and decorrelation time, allowing the model to have specified properties or for the model to be tuned to data, either to real data or a free run of sophisticated ocean model.

\begin{figure}
    \begin{center}
        \includegraphics[]{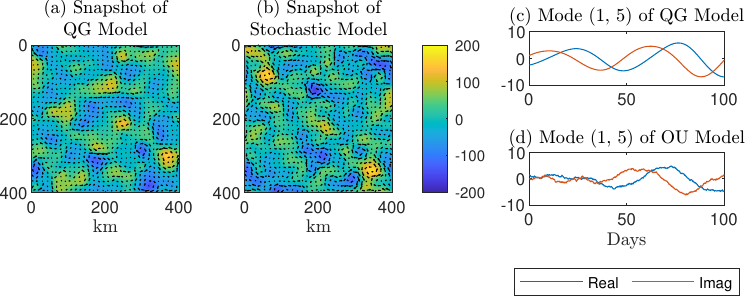}
    \end{center}
    \caption{\label{fig:qgvsou}Qualitative comparison between the simulations from the quasi-Geostrophic ocean model and from the stochastic ocean model. The simulations for each model are independent of each other. Therefore, the comparison is not based on the point-wise error between the two simulation but instead on the qualitative similarity. Panels (a) and (b): snapshots of the stream function and associated velocity field of the two models, respectively. Panels (c) and (d): the real and imaginary parts of mode $(1, 5)$ of the two models.
    }
\end{figure}

Figure \ref{fig:qgvsou} shows a qualitative comparison between the simulations from the quasi-Geostrophic ocean model and from the stochastic ocean model. The simulations for each model are independent of each other. Therefore, the comparison is not based on the point-wise error between the two simulation but instead on the qualitative similarity.

\section{An efficient analytically solvable nonlinear Lagrangian data assimilation framework\label{sec:cg}}
The DEM ice floe model and the two-layer QG ocean model are utilized to generate the true signal of both the floe trajectories and the ocean field. On the other hand, to facilitate an efficient data assimilation scheme, the stochastic models are coupled with the DEM ice floe model as the forecast model. Furthermore, using an approximation where $\mathbf{F}_\mathrm{o}^\ell$ and $\tau_\mathrm{o}^\ell$ are uniform over the floe based on the value at the center of mass, then the DEM ice floe model can have a more explicit form,
\begin{equation}\label{eqn:floe_simplified}
\begin{split}
    \frac{\d\mathbf{x}^\ell}{\d t} =& \mathbf{u}^\ell \\
      \frac{\d\mathbf{u}^\ell}{\d t} =&  \frac{d_\mathrm{o} \rho_\mathrm{o}}{h_\ell \rho_{ice}} \left(\mathbf{u}_\mathrm{o}(\mathbf{x}^\ell, t) - \mathbf{u}^\ell\right)\left|\mathbf{u}_\mathrm{o}(\mathbf{x}^\ell, t) - \mathbf{u}^\ell\right| \\
    \frac{\d\Omega^\ell}{\d t} =& \omega^\ell \\
     \frac{\d\omega^\ell}{\d t} =&  \frac{d_\mathrm{o} \rho_\mathrm{o}}{h_\ell \rho_\mathrm{ice}} \left(\frac{\nabla \times \mathbf{u}_\mathrm{o}(\mathbf{x}^\ell, t)}{2} - \omega\right)\left|\frac{\nabla \times \mathbf{u}_\mathrm{o}(\mathbf{x}^\ell, t)}{2} - \omega\right|. \\
\end{split}
\end{equation}
Such a setup has been utilized in the previous work \cite{damsgaard2018application, chen2021lagrangian} when the ice floes are assumed to have cylindrical shapes with uniform thickness. Nevertheless, it can also be applied to the case in this work where floes have arbitrary shapes as long as the recovery of the ocean focuses on spatial scales larger than the floe size, which is the case for the mesoscale eddies here. Note that this further simplification is not necessary to develop an analytically solvable data assimilation framework, but it accelerates the computations.

Aggregating the DEM floe model \eqref{eqn:floe_simplified} with the linear stochastic ocean forecast model \eqref{eqn:ocean_velocity}, \eqref{eqn:ocean_vorticity} and \eqref{eqn:ouprocess} leads to the coupled forecast system for Lagrangian data assimilation. After approximating the quadratic drag terms by a linear drag with time-varying coefficients, the coupled forecast system remains highly nonlinear due to displacement and angular displacement variables appearing in the exponential functions. Nevertheless, data assimilation uses the observations of these variables to estimate the velocity of the floes and the ocean. Conditioned on the observed displacement and angular displacement of each floe (meaning they are given), the equations of the ocean and floe velocity fields are conditionally linear, and the associated distribution of these variables is conditional Gaussian. For such a system, closed analytic formulae are available to compute the conditional distribution (also known as the posterior distribution) of these state variables. The mathematical details are included as follows.

A coupled system falls under such a framework is called a conditional Gaussian (CG) nonlinear system. It can be written in the following form \cite{liptser2013statistics, chen2022conditional}:
\begin{align}
    \label{eqn:observed} \frac{\d\mathbf{X}}{\d t} =& \mathbf{A}_0(\mathbf{X}, t) + \mathbf{A}_1(\mathbf{X}, t) \mathbf{Y} + \mathbf{B}(\mathbf{X}, t) \dot{\mathbf{W}}_\mathbf{X}(t), \\
    \label{eqn:unobserved} \frac{\d\mathbf{Y}}{\d t} =& \mathbf{a}_0(\mathbf{X}, t) + \mathbf{a}_1(\mathbf{X}, t) \mathbf{Y} + \mathbf{b}(\mathbf{X}, t) \dot{\mathbf{W}}_\mathbf{Y}(t),
\end{align}
where $\mathbf{X}$ is the state vector of observed variables and $\mathbf{Y}$ is the state vector of unobserved variables.
Here, $\mathbf{A}_0$, $\mathbf{A}_1$, $\mathbf{B}$, $\mathbf{a}_0$, $\mathbf{a}_1$, and $\mathbf{B}$ are matrices that can depend on $t$ and $\mathbf{X}$ arbitrarily nonlinearly.
While this dependence is always assumed, it may not be denoted in the subsequent discussion for notational efficiency.
$\dot{\mathbf{W}}_\mathbf{X}$ and $\dot{\mathbf{W}}_\mathbf{Y}$ are Gaussian white noise. In the coupled floe-ocean system, $\mathbf{X}$ contains the observed displacement and angular displacement of each individual floe. The state variable $\mathbf{Y}$ is a collection of the velocities and angular velocities of the floes as well as the Fourier coefficients of the ocean fields. Because of the potentially highly nonlinear interactions with the observed variables, the CG framework has been applied to a wide range of nonlinear systems and an even wider range of systems have suitable nonlinear CG approximations \cite{chen2018conditional}.

Despite the nonlinearity, all the matrices ($\mathbf{A}_1$, $\mathbf{a}_1$, $\mathbf{B}$ and $\mathbf{b}$) and vectors ($\mathbf{A}_0$ and $\mathbf{a}_0$) only depend on the observed variables, $\mathbf{X}$, and do not depend on the unobserved variables, $\mathbf{Y}$. Therefore, the system is conditionally linear in $\mathbf{Y}$ given  $\mathbf{X}$. This property ensures that the conditional distribution of $\mathbf{Y}$ given a trajectory of $\mathbf{X}$ is Gaussian, as the name implies. Because of this property, the posterior distributions for the data assimilation solution, including both the filtering and the smoothing, are characterized by their posterior mean vectors and covariance matrices, which are given by closed analytic formulae.

The posterior distribution of $\mathbf{Y}(t)$ at time $t$ given a trajectory of $\mathbf{X}(s)$ over the interval $[0, t]$ is known as the filter posterior distribution and is given by
\begin{equation}
    p(\mathbf{Y}(t) \mid \mathbf{X}(s), s \leq t) \sim \mathcal{N}(\boldsymbol{\mu}_\mathrm{f}(t), \mathbf{R}_\mathrm{f}(t))
\end{equation}
where $\mathbf{\boldsymbol{\mu}_\mathrm{f}}$ and $\mathbf{R}_\mathrm{f}$ are the filter mean and covariance respectively.
The filter mean and covariance can be calculated using the forward equations
\begin{align}
    \label{eqn:filtermean} \frac{\d\boldsymbol{\mu}_\mathrm{f}}{\d t} =& (\mathbf{a}_0 + \mathbf{a}_1 \boldsymbol{\mu}_\mathrm{f}) + (\mathbf{R}_\mathrm{f} \mathbf{A}_1^\ast) (\mathbf{B} \mathbf{B}^\ast)^{-1} \left(\frac{\d\mathbf{X}}{\d t} - (\mathbf{A}_0 + \mathbf{A}_1 \boldsymbol{\mu}_\mathrm{f})\right) \\
    \label{eqn:filtercov} \frac{\d\mathbf{R}_\mathrm{f}}{\d t} =& \mathbf{a}_1 \mathbf{R}_\mathrm{f} + \mathbf{R}_\mathrm{f} \mathbf{a}_1^\ast + \mathbf{b}\mathbf{b}^\ast - (\mathbf{R}_\mathrm{f} \mathbf{A}_1^\ast)(\mathbf{B} \mathbf{B}^\ast)^{-1} (\mathbf{A}_1 \mathbf{R}_\mathrm{f})
\end{align}
where the initial condition is a Gaussian distribution with mean $\boldsymbol{\mu}_\mathrm{f}(0)$ and covariance $\mathbf{R}_\mathrm{f}(0)$.
Here ``$\cdot^\ast$'' denotes the conjugate transpose.

The posterior distribution of $\mathbf{Y}(t)$ at time $t \in [0, T]$ given a trajectory of $\mathbf{X}(s)$ over the entire interval $[0, T]$ is known as the smoother posterior distribution \cite{chen2020efficient} and is given by
\begin{equation}
    p(\mathbf{Y}(t) \mid \mathbf{X}(s), s \in [0, T]) \sim \mathcal{N}(\boldsymbol{\mu}_\mathrm{s}(t), \mathbf{R}_\mathrm{s}(t))
\end{equation}
where $\mathbf{\boldsymbol{\mu}_\mathrm{s}}$ and $\mathbf{R}_\mathrm{s}$ are the smoother mean and covariance respectively.
Compared to the filter posterior distribution, the smoother posterior distribution incorporates both past and future observational observation.

To calculate the smoother mean and covariance, first the filter mean and covariance are calculated.
The condition of the smoother mean and covariance equations at the final time $T$ is given by the filter mean and covariance at time $T$, that is $(\mathbf{\boldsymbol{\mu}_\mathrm{s}}(t), \mathbf{R}_\mathrm{s}(t)) = (\mathbf{\boldsymbol{\mu}_\mathrm{f}}(t), \mathbf{R}_\mathrm{f}(t))$. The following equations are then integrated backwards in time from time $T$ to time $0$:
\begin{align}
    \label{eqn:smoothermean} \frac{\d \boldsymbol{\mu}_\mathrm{s}}{\d t} =& -\mathbf{a}_0 -\mathbf{a}_1 \boldsymbol{\mu}_\mathrm{s} + (\mathbf{b} \mathbf{b}^\ast) \mathbf{R}_\mathrm{f}^{-1} (\boldsymbol{\mu}_\mathrm{f} - \boldsymbol{\mu}_\mathrm{s})\\
    \label{eqn:smoothercov} \frac{\d \mathbf{R}_\mathrm{s}}{\d t} = & -(\mathbf{a}_1 + (\mathbf{b} \mathbf{b}^\ast) \mathbf{R}_\mathrm{f}^{-1}) \mathbf{R}_\mathrm{s} - \mathbf{R}_\mathrm{s} (\mathbf{a}_1^\ast + (\mathbf{b} \mathbf{b}^\ast) \mathbf{R}_\mathrm{f}^{-1}) + \mathbf{b} \mathbf{b}^\ast.
\end{align}
The smoother equations calculate the matrix inversion $\mathbf{R}_\mathrm{f}^{-1}$ which can be computationally inconvenient for high dimensional systems.

While the smoother mean and covariance can be used to sample $\mathbf{Y}$ at a fixed time instant, the CG framework can also be used to sample entire trajectories of $\mathbf{Y}$ conditioned on $\mathbf{X}$ on the interval $[0, T]$. Such sampled trajectories are essential for calculating eddy statistics, such as the eddy lifetime, or other eddy diagnostics, such as the Lagrangian descriptor, incorporating temporal information from the velocity field. Sampling these trajectories is accomplished via the backward sampling equation where an initial $\mathbf{Y}(T)$ is drawn from $\mathbf{Y}(T) \sim \mathcal{N}(\boldsymbol{\mu}_\mathrm{f}(T), \mathbf{R}_\mathrm{f}(T))$ and then its trajectory is calculated using the following stochastic equation integrated backward in time
\begin{equation}
    \label{eqn:backwardsampling} \frac{\d \mathbf{Y}}{\d t} = -\mathbf{a}_0 -\mathbf{a}_1 \mathbf{Y} + (\mathbf{b} \mathbf{b}^\ast) \mathbf{R}_\mathrm{f}^{-1} (\boldsymbol{\mu}_\mathrm{f} - \mathbf{Y}) + \mathbf{b} \dot{\mathbf{W}}_\mathbf{Y}(t).
\end{equation}
When sampling multiple trajectories of $\mathbf{Y}$ conditioned on $\mathbf{X}$, the same filter mean and covariance are used in equation \eqref{eqn:backwardsampling} and the dynamics of the sampled trajectories differ only in the sampled initial condition at time $T$ and the realizations of the Gaussian white noise $\dot{\mathbf{W}}_\mathbf{Y}(t)$. In particular, $\mathbf{R}_\mathrm{f}^{-1}$ can be calculated once on the interval $[0, T]$ and then reused to calculate each sampled trajectory. Note that various additional techniques, such as reduced-order data assimilation schemes and online smoothers building upon the above framework, help further reduce the computational cost and storage for systems with more complexity and higher dimensionality \cite{chen2023uncertainty}.

\section{\label{sec:evofow} Expected value of the OW parameter}

The OW parameter, as most eddy diagnostics, is a nonlinear function of the ocean velocity.
Because of this nonlinearity, calculating the OW parameter from the mean estimate of the ocean state is different from the expected value of the OW parameter under the entire posterior distribution. The contribution of the uncertainty present in the posterior distribution can be clearly illustrated by calculating the expected value of the OW parameter explicitly. Note that, for notation simplicity, the two components of $\mathbf{u}_o$, namely $(u_o,v_o)$, are written as $(u,v)$ hereafter, which should not be confused with the floe velocity in \ref{sec:ice}.

The OW parameter is given by
\begin{align}
    \mathrm{OW}(\mathbf{u}_o) =& (u_x - v_y)^2 + (v_x + u_y)^2 - (v_x - u_y)^2 \\
                =& (u_x - v_y)^2 + 4v_x u_y \\
                =& u_x^2 - 2u_xv_y + v_y^2 + 4v_x u_y,
\end{align}
where a subscript indicates a partial derivative in space.
To calculate the expected value of the OW parameter we apply a mean fluctuation decomposition
\begin{equation}
    u_x = \bar{u}_x + u_x',\qquad u_y = \bar{u}_y + u_y',\qquad  v_x = \bar{v}_x + v_x',\qquad\mbox{and}\qquad v_y = \bar{v}_y + v_y',
\end{equation}
where $\bar{u}_x$, $\bar{u}_y$, $\bar{v}_x$, and $\bar{v}_y$ are the mean of each variable and $u_x'$, $u_y'$, $v_x'$, and $v_y'$ are random variables with mean 0.
Then the OW parameter can be written as
\begin{align}
    \mathrm{OW}(\mathbf{u}_o) =& (\bar{u}_x + u_x'  - (\bar{v}_y + v_y'))^2 + 4(\bar{v}_x + v_x') (\bar{u}_y + u_y') \\
    =& \bar{u}_x\bar{u}_x + \bar{u}_x u_x'  - \bar{u}_x \bar{v}_y - \bar{u}_x v_y' + u_x'\bar{u}_x + u_x'u_x'  - u_x'\bar{v}_y - u_x'v_y' \\
    &  -\bar{v}_y\bar{u}_x -\bar{v}_yu_x' + \bar{v}_y\bar{v}_y +\bar{v}_y v_y' - v_y'\bar{u}_x - v_y'u_x'  + v_y'\bar{v}_y + v_y'v_y'\\
    & + 4(\bar{v}_x\bar{u}_y + \bar{v}_xu_y' + v_x'\bar{u}_y + v_x'u_y').
\end{align}
Calculating the expected value yields
\begin{align}
    \mathbb{E}[\mathrm{OW}(\mathbf{u}_o)] =& \bar{u}_x^2  - 2\bar{u}_x \bar{v}_y + \bar{v}_y^2 + 4\bar{v}_x\bar{u}_y \\
    &  + \mathbb{E}[(u_x')^2] - 2\mathbb{E}[u_x'v_y'] + \mathbb{E}[(v_y')^2] + 4\mathbb{E}[v_x'u_y'] \\
    =& \mathrm{OW}(\bar{\mathbf{u}}_o) + \mathbb{E}[(u_x')^2] - 2\mathbb{E}[u_x'v_y'] + \mathbb{E}[(v_y')^2] + 4\mathbb{E}[v_x'u_y'].
\end{align}
Then the expected value of the OW parameter is the OW parameter of the posterior mean with additional terms for the contribution of the uncertainty.

The above calculation works for a general ocean model.
In the case that the ocean is incompressible, then $u_x = -v_y$ and so the above simplifies to
\begin{equation}
    \mathrm{OW}(\mathbf{u}_o) = 4u_x^2 + 4v_x u_y
\end{equation}
and
\begin{equation}
    \mathbb{E}[\mathrm{OW}(\mathbf{u}_o)] = \mathbb{E}[\mathrm{OW}(\bar{\mathbf{u}}_o)] + 4\mathbb{E}[(u_x')^2] + 4\mathbb{E}[v_x'u_y'].
\end{equation}

\section{\label{sec:twinexperiment} Eddy identification in a perfect model twin experiment}
This section includes the eddy identification results in a perfect model twin experiment. The model generating the true signal and the model used for state estimation are both linear stochastic models. This contrasts with the setup in the results shown in the main text, where the underlying truth is given by the two-layer QG model while the linear stochastic model is used as an approximate forecast model for state estimation ocean. Consistent with the setup in the main text, the setup here contains the use of Fourier modes $\mathbf{k}\in[-11,11]^2$, which includes $23^2 - 1 = 528$ modes (the $(0, 0)$ mode is identically zero). The parameters in the linear stochastic modes are the same as the calibrated ones based on the QG model. The purpose of such a perfect model experiment is to understand the model error affecting the probabilistic eddy identification forecast.

\begin{figure}
    \begin{center}
        \includegraphics[]{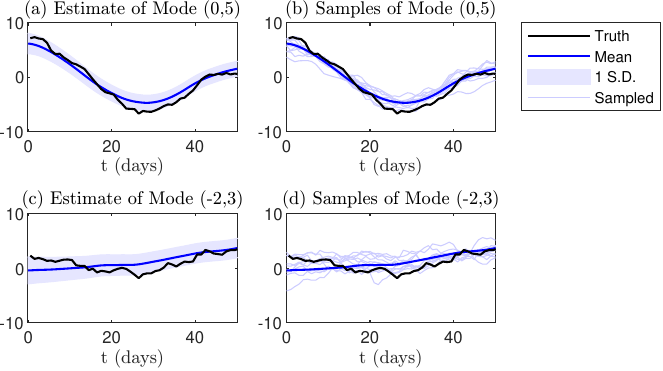}
    \end{center}
    \caption{\label{fig:smoother_twin_experiment}Eddy identification in a twin experiment setup. The caption is similar to Figure \ref{fig:smoother}.
    }
\end{figure}

\begin{figure}
    \begin{center}
        \includegraphics[]{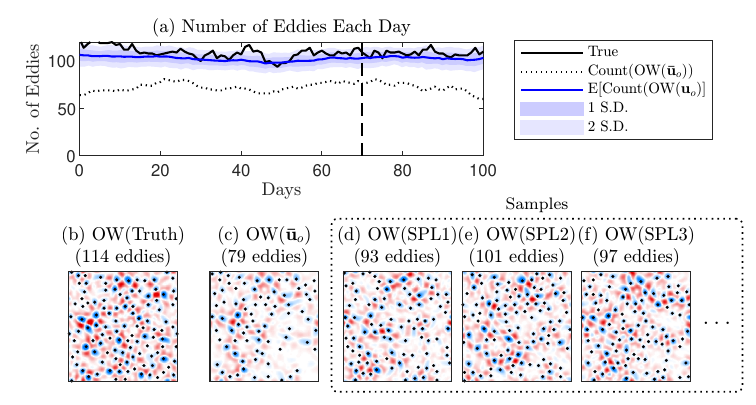}
    \end{center}
    \caption{\label{fig:count_fig_twin_experiment}Eddy identification in a twin experiment setup. The caption is similar to Figure \ref{fig:count}.
    }
\end{figure}

\begin{figure}
    \begin{center}
        \includegraphics[]{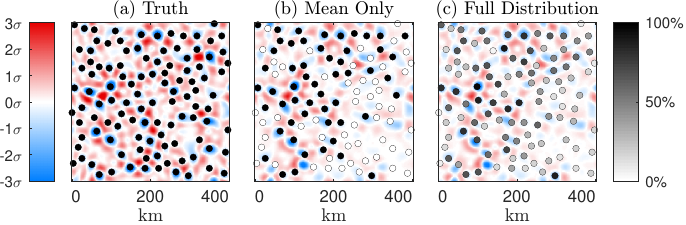}
    \end{center}
    \caption{\label{fig:scatter_fig_twin_experiment}Eddy identification in a twin experiment setup. The caption is similar to Figure \ref{fig:identification}.
    }
\end{figure}

Figures \ref{fig:smoother_twin_experiment}--\ref{fig:scatter_fig_twin_experiment} are the analogs to Figures \ref{fig:smoother}, \ref{fig:count}, and \ref{fig:identification} in the main text. The results in these figures are qualitatively similar to those in the main text, justifying the appropriateness and robustness of the framework in the presence of small model errors. The only noticeable difference is between Figure \ref{fig:count} and Figure \ref{fig:count_fig_twin_experiment}. The mean of the PDF of the eddy number, $\mathbb{E}[\mbox{Count}((\mbox{OW}({\mathbf{u}}_o))]$ (blue solid curve), is overestimated in the former case. In contrast, it is quite consistent with the truth in the latter. As a result, there are very few eddies in Panels (d)--(f) of Figure \ref{fig:count_fig_twin_experiment} that are falsely identified, which differs from the case in Figure \ref{fig:count}. The main reason for the difference between the two experiments is the following: Although the linear stochastic model is calibrated based on the forecast statistics of the QG model, the noise induces slightly more temporal randomness, which can be seen in Panels (c)--(d) of Figure \ref{fig:qgvsou}. As the OW parameter is based on individual snapshots, such randomness causes more variability and additional eddies.

%
%

\section*{Open Research Section}
The code used to process the data and create the figures was written in MATLAB. The code and output data of the experiments are available on Zenodo: \\https://doi.org/10.5281/zenodo.11188124

\section*{Acknowledgement}
The research of N.C. is funded by Naval Research (ONR) Multidisciplinary University Research Initiative (MURI) award N00014-19-1-2421 and ONR N00014-24-1-2244. J.C. is supported as a research assistant/research scientist under these grants. S.W. acknowledges the financial support provided by the EPSRC Grant No. EP/P021123/1 and the support of the William R. Davis '68 Chair in the Department of Mathematics at the United States Naval Academy. The research of E.L. is supported by ONR N0001423WX01622.

\bibliographystyle{elsarticle-num}

\end{document}